\documentclass{elsart}
\usepackage[latin1]{inputenc}
\usepackage{latexsym}
\usepackage{graphicx}
\usepackage{amsfonts}
\usepackage{amsmath}
\usepackage{amssymb}
\theoremstyle{plain}
\newtheorem{nota}[thm]{Remark}
\newtheorem{ejem}[thm]{Example}

\newcommand{\N}{\mathbb{N}}

\newcommand{\R}{\mathbb{R}}

\newcommand{\luno}{\ell^1}
\newcommand{\ltwo}{\ell^2}

\newcommand{\f}{\varphi}

\newcommand{\ASP}{{ASP}}
\newcommand{\eps}{\varepsilon}
\newcommand{\limty}[1]{\lim\limits_{#1\to\infty}}
\newcommand{\limn}{\limty{n}}
\newcommand{\bx}{\hfill$\blacksquare$}
\newcommand{\wbx}{\hfill$\square$}

\newcommand{\dem}{\noindent{\bf Proof. }}

\def\mles{{$M$-LES}}
\def\ces{{CES}}
\def\sces{SCES}

\begin{document}

\begin{frontmatter}
\title{Asymptotic almost-equivalence of abstract evolution systems\thanksref{funds}}
\thanks[funds]{Supported by the Millennium Institute on Complex Engineering Systems (MIDEPLAN P05-004-F),
FONDECYT 1050706, MECESUP UCH0009 and Basal Proyect CMM, Universidad de Chile.}

\author{Felipe Alvarez}
\address{CMM (CNRS UMI 2807), Departamento de Ingenier\'\i a Matem\'atica, Universidad de Chile, Av. Blanco Encalada 2120, Santiago, Chile.}
\ead{falvarez@dim.uchile.cl}

\author{Juan Peypouquet\corauthref{cor}}
\corauth[cor]{Corresponding author.}
\address{Departamento de Matem\'atica, Universidad T\'ecnica Federico Santa Mar\'\i a, Av. Espa\~na 1680, Valpara\'\i so, Chile.}
\ead{juan.peypouquet@usm.cl}

\begin{abstract}
We study the asymptotic behavior of almost-orbits of abstract evolution systems in Banach spaces with or without a
Lipschitz assumption. In particular, we establish convergence, convergence in average and
almost-convergence of almost-orbits both for the weak and the strong topologies based on the behavior of the orbits. We
also analyze the set of almost-stationary points.

\end{abstract}


\end{frontmatter}

\section{Introduction and preliminaries}

Roughly speaking, a dynamical system in discrete (resp. continuous) time is a rule that
determines a sequence (resp. trajectory) departing from certain initial data and which
evolves in an either finite or infinite dimensional space. In this sense, any iterative
algorithm  may be con\-si\-dered as a discrete evolution system.  If it is possible to
find a continuous-in-time version for the discrete procedure, it is then natural to
expect that some of the properties of the former are close to the similar properties of
the latter. Of course, in general such an inheritance of properties is not true without
additional conditions, in particular on the parameters of the algorithm. The dynamical
approach to iterative methods has certain advantages: a continuous-in-time evolution
system satisfying nice qualitative properties may suggest new iterative methods, and
sometimes the techniques used to investigate the continuous case can be adapted to obtain
results for the discrete algorithm. On the other hand, one may be concerned with
different aspects of the trajectories of a given continuous-in-time dynamical system,
namely: existence, exact or approximate computation, regularity, long-term behavior,
stability, numerical integration. Of special interest is the effect of certain
perturbations of the original system on the qualitative properties of the corresponding
trajectories. In this context, this paper deals with some of the asymptotic properties
that are common to systems which can be considered equivalent in a sense  to be made
precise later on.

Let $C$ be a nonempty Borel subset of a Banach space $(X, \|\cdot\|)$. An {\em evolution system} (ES for short) on $C$
is a two-parameter family $U=\{U(t,s) \mid  t\ge s\ge 0\}$ of possibly non-linear maps from $C$ into itself satisfying:

\begin{itemize}
        \item [i)] $\forall t\geq 0$, $\forall x\in C$, $U(t,t)x=x$; and
        \item [ii)] $\forall t\geq s\geq r\geq 0$, $\forall x\in C$, $U(t,s)U(s,r)x=U(t,r)x$.
\end{itemize}
A {\em $M$-Lipschitz evolution system} (\mles) is an ES $U$ such that
$\|U(t,s)x-U(t,s)y\|\le M\|x-y\|$ for some $M>0$ and all $t\geq s\geq  0$, $x,y\in C$. A
{\em contracting evolution system} (\ces) is a 1-LES.

An ES $U$ is {\em autonomous} if for all $t,s\geq 0$ we have $U(t,0)=U(t+s,s)$. For such
an ES, the family  $T=\{ T(t):=U(t,0)\mid t\geq 0\}$ defines a {\em semigroup}, that is
$T(0)x=x$ and $T(t)T(s)x=T(t+s)x$ for all $t,s \geq 0$, $x\in C$.

\begin{ejem}\label{E:es1}{\em Let $F$ be a (possibly multivalued) function from $[t_0,\infty)\times
C$ to $C$. Suppose that for every $s\ge t_0$ and $x\in C$ the differential inclusion
$u'(t) \in F(t,u(t))$, with initial condition $u(s) = x$, has a unique solution
$u_{s,x}:[s,\infty)\mapsto C$. The family $U$ defined by $U(t,s)x=u_{s,x}(t)$ is an
evolution system on $C$. If $X$ is Hilbert space and $F(t,x)=-A_tx$, where $\{A_t\}$ is a
family of maximal monotone operators, then the corresponding $U$ is a \ces. }\wbx
\end{ejem}

\begin{ejem}{\em \label{E:es2} Take a strictly increasing unbounded sequence $\{\sigma_n\}$ of
positive numbers and set $\nu(t)=\max\{n\in\N\ |\ \sigma_n\le t\}$. Consider a family
$\{F_n\}$ of functions from $C$ into $C$ and define
$U(t,s)=\prod_{n=\nu(s)+1}^{\nu(t)}F_n$, the product representing composition of
functions. Then $U$ is an ES. If each $F_n$ is $M_n$-Lipschitz and the product
$\prod_{n=1}^\infty M_n$ is bounded from above by $M$, then $U$ is an \mles. For
instance, if $F_n=(I+A_n)^{-1}$, where $\{A_n\}$ is a family of $m$-accretive operators
on $C$, then the piecewise constant interpolation of infinite products of resolvents
defines a \ces.} \wbx\end{ejem}

If $U$ is an ES on $C$, an {\em orbit} of $U$ is a function  $u:[0,\infty)\to C$  such that for some  $t_0\geq 0$ and
$x_0\in C$, $u(t)=U(t,t_0)x_0$ for all $t\geq t_0$.  Throughout this paper, all orbits are assumed to be measurable and
locally bounded, hence locally integrable on $[0,\infty)$. More generally, we say that a function $u\in
L^\infty_{\hbox{\tiny loc}}(0,\infty;C)$ is an {\em almost-orbit} of $U$ if
\vskip-20pt
\begin{equation}\label{ao} \limty{t}\sup_{h\geq 0}\|u(t+h)-U(t+h,t)u(t)\|=0.
\end{equation}

Clearly, orbits are almost-orbits. Two ES are {\em asymptotically almost-equivalent}
(AAE) if every orbit of each one is an almost-orbit of the other.

If $U$ is an autonomous \ces\ and $V$ is AAE to $U$, then $V$ is an {\em asymptotic
semigroup} as defined in \cite{Pas}. In that case, every orbit of $U$ converges strongly
(or weakly) if, and only if, every orbit of $V$ does (see Proposition \ref{convmles}
below).

\begin{nota}\label{R:limits}{\em Suppose $U$ is a \mles\ on $C$. If $u$ is an almost-orbit of $U$, then so is any
function $v\in L^\infty_{\hbox{\tiny loc}}(0,\infty;C)$ satisfying
$\limty{t}\|v(t)-u(t)\|=0$. }\wbx\end{nota}

\begin{nota}{\em Suppose that for each $r>0$ there is $G_r:\R_+^2\to\R_+$ such that
$\limty{t}\sup\limits_{h\ge 0}G_r(t+h,t)=0$ and $\|U(t,s)x-V(t,s)x\|\le G_r(t,s)$ for all
$x\in B(0,r)$ and $t\ge s\ge 0$. Every bounded orbit of $U$ is an almost-orbit of $V$ and
viceversa. If the same $G_r\equiv G$, the boundedness assumption is
unnecessary.}\wbx.\end{nota}

The term ``almost-orbit" was introduced in \cite{KoM} for a {\em continuous} function
satisfying \eqref{ao}. Later, in \cite{Miy}, the author gives a weaker definition, just
requiring $\limty{t,h}\|u(t+h)-U(t+h,t)u(t)\|=0,$ but still for continuous functions. The
latter is slightly weaker than \eqref{ao} for practical purposes. In fact, the example
provided in \cite{Miy} to motivate the interest of studying almost-orbits also satisfies
\eqref{ao}. In both cited works the authors give criteria that can be applied to an
almost-orbit in order to ensure certain asymptotic behavior. The same approach is used in
\cite{Xu}. In \cite{KiR,Rou}, the authors carry out a similar analysis for {\em uniformly
asymptotically almost nonexpansive curves} (a concept that includes almost-orbits {\em of
almost nonexpansive semigroups}) in Hilbert space. Other results on the asymptotic
behavior of almost-orbits of nonexpansive semigroups can be found in \cite{LNT} (see also
the references therein). Notice that \cite{KoM,Miy} contain versions of Proposition
\ref{almostmles} below in different settings. Our intention is to show how to derive many
asymptotic properties of almost-orbits by studying only the orbits. In that sense our
work is different but complementary to \cite{KoM,Miy}.

\begin{ejem}{\em Let $A$ be a $m$-accretive operator on $X$ and let $U$ be the
autonomous \ces\ defined by the inclusion $-\dot u\in Au$ as in Example \ref{E:es1}.}
\begin{enumerate}
    \item {\em Take $f\in L^1(0,\infty;X)$ and let $V$ be the \ces\ defined by the integral solutions of $-\dot u\in Au+f$.
        The orbits of $V$ are almost-orbits of $U$ (see \cite{KoM}). This result is generalized in
        \cite{AlP1}.}
    \item {\em For a sequence $\{\lambda_n\}$ in $(0,\infty)$ set $\sigma_n=\sum_{k=0}^n\lambda_k$ and $\nu(t)$ as
        in Example \ref{E:es2}. Define the \ces\ $W(t,s)=\prod_{n=\nu(s)+1}^{\nu(t)}(I+\lambda_nA)^{-1}$. If
        $\{\lambda_n\}\in\ltwo\setminus\luno$ then $U$ and $W$ are AAE; see \cite{KoS},
        although the fact that the orbits of $W$ are almost-orbits of $U$ had already been proved in \cite{KoM}.
        This was shown earlier in \cite{Pas} by  assuming $A$ to be single-valued and
        Lipschitz.}
    \item {\em If $A$ is the subdifferential of a proper, closed and convex function in a Hilbert
        space, $U$ and $W$ are AAE if $\{\lambda_n\}\notin\luno$ (see \cite{Gul})}.\wbx
\end{enumerate}
\end{ejem}

We shall see that orbits and almost-orbits of an ES have the same asymptotic behavior in
terms of boundedness, convergence and other related properties. A few results of this
kind can be found in \cite{Pas,KoM,KoS,Gul} for differential inclusions of the type
$-u'(t)\in Au(t)$, where $A$ is $m$-accretive. A first attempt to deal with nonautonomous
and non-Lipschitz systems can be found in \cite{Miy}.

The paper is organized as follows: In section \ref{S:mles} we focus on \mles. We show
basic properties of their almost-orbits, exploiting the dissipativity behind the uniform
Lipschitz constant. We also state and prove some asymptotic equivalence results. Despite
the surprising fact that the main convergence results hold for arbitrary ES almost as
stated, we prefer to present this easier case first for the sake of a clear exposition.
Section 3 contains further results on the special class of strongly contracting ES and
the structure of the set of almost-stationary points of an \mles. The asymptotic
properties for general ES without any assumptions on the space-dependence are given in
section 4. Neither Lipschitz continuity nor asymptotic nonexpansiveness is imposed. We
present some additional results on uniform continuity and cluster points in section 5 and
some remarks on the applicability of this theory in section 6.

\section{Lipschitz evolution systems}\label{S:mles}

The following is similar to \cite[Lemma 3.1]{KoM}:
\begin{prop}\label{prop1}
Let $U$ be a \mles\ and $u_1,u_2$ two almost-orbits of $U$. Then
\begin{itemize}
    \item [i)] $\limsup\limits_{t\to\infty}\|u_1(t)-u_2(t)\|\le M\liminf\limits_{t\to\infty}\|u_1(t)-u_2(t)\|<\infty.$
    \item [ii)] If one almost-orbit of $U$ is bounded, then every almost-orbit of $U$ is.
    \item [iii)] If 0 is a cluster point of $\|u_1(t)-u_2(t)\|$ then $\limty{t}\|u_1(t)-u_2(t)\|=0$.
    \item [iv)] If $U$ is a \ces, the limit $\limty{t}\|u_1(t)-u_2(t)\|$ always exists.
\end{itemize}
\end{prop}
\dem We just prove $i)$. For $i=1,2$ let $\psi_i(t)=\sup\limits_{h\ge
0}\|u_i(t+h)-U(t+h,t)u_i(t)\|$.
Then $\|u_1(t+h)-u_2(t+h)\| \le \psi_1(t)+\psi_2(t)+M\|u_1(t)-u_2(t)\|$ for every $h\geq
0$. Hence $\limsup\limits_{h\to\infty}\|u_1(h)-u_2(h)\|\le
\psi_1(t)+\psi_2(t)+M\|u_1(t)-u_2(t)\|<\infty$ and finally
$\limsup\limits_{h\to\infty}\|u_1(h)-u_2(h)\|\le
M\liminf\limits_{t\to\infty}\|u_1(t)-u_2(t)\|.$ \bx

\begin{nota}\label{R:all_bounded}{\em Let $U$, $V$ be \mles\ which are AAE. If {\bf one}
almost-orbit of $U$ {\bf or} $V$ is bounded, {\bf every} almost-orbit of $U$ {\bf and}
$V$ is so. }\wbx\end{nota}

%
\begin{nota}\label{boundbound}{\em Let $u$ be a bounded almost-orbit of an \mles\ $U$ so that
$\|u\|_\infty=\sup_t \|u(t)\|<\infty$. Since $u$ is an almost-orbit, there exists
$p_0\geq 0$ such that for all $p\geq p_0$, we have $\|u(p+h)-U(p+h,p)u(p)\|\leq 1\hbox{
for all } h\geq 0.$ Hence, for all $p\geq p_0$ and $h\geq 0$ we get $\|U(p+h,p)u(p)\|\leq
1+\|u\|_\infty$. }\wbx\end{nota}


The following result and its proof are inspired by \cite[Lemma 1]{Pas}, where the author studies two special cases: when $U$ is an autonomous \ces, and when the almost-orbits are in fact the orbits of a semigroup of contractions. This result had already been presented by the authors for arbitrary \ces\ in \cite{AlP1}. Here we give a shorter proof in a more general context.

\begin{prop}\label{convmles}
Let $U$ be an \mles. If every orbit of $U$ converges strongly (weakly) as time goes to
infinity, so does every almost-orbit of $U$.
\end{prop}
\dem Let $\tau$ denote the hypothesized topology. And suppose that
the $\tau-$limit of $U(t,s)x$ as $t\to\infty$ exists for all $x$ and
$s$. Let $u$ be an almost-orbit of $U$. Take $p\ge 0$ and set
$\zeta(p)=\tau-\limty{t}U(t,p)u(p)$. We have\vskip-20pt
$$\zeta(p+h)-\zeta(p)=\tau-\limty{t}\{U(t,p+h)u(p+h)-U(t,p)u(p)\}.$$
But for all $t\geq p+h$ the quantity
$\| U(t,p+h)u(p+h)-U(t,p)u(p)\|$ can be bounded above by $M\| u(p+h)-U(p+h,p)u(p)\|$
and by $\tau$-lower semicontinuity of the norm we get $\|\zeta(p+h)-\zeta(p)\|\le M\| u(p+h)-U(p+h,p)u(p)\|.$
Since $u$ is an almost-orbit of $U$, the right-hand side  tends to zero as $p\to\infty$ uniformly in $h\ge 0$. Therefore $\{\zeta(p):p\to \infty\}$ is a Cauchy net that converges strongly
to a limit $\zeta_\infty$. Finally, we can express $u(p+h)-\zeta_\infty$, for all $p,h\geq 0$, as
$u(p+h)-\zeta_\infty  = [u(p+h)-U(p+h,p)u(p)]+[U(p+h,p)u(p)-\zeta(p)]+[\zeta(p)-\zeta_\infty].$
Given $\eps>0$ we can choose $p$ large enough so that the first and third terms on the right-hand side are less than $\eps$ in norm, uniformly in $h$ for the first term. Next for such a fixed $p$, we let $h\to\infty$ so that the second term $\tau-$converges to zero.  Hence $u(t)$ is $\tau-$convergent to $\zeta_\infty$ as $t\to\infty$. \bx

\begin{nota}{\em Consider the following more general setting: Let $(X,d)$ be a complete metric space (not even the linear structure is necessary). The Lipschitz condition in the definition of \mles\ reads $d(U(t,s)x,U(t,s)y)\le Md(x,y)$. The definition
of almost-orbit can be rephrased as $\limty{t}\sup_{h\geq 0}d(u(t+h),U(t+h,t)u(t))=0$. It
is easy to see that Proposition \ref{prop1} and the statement in Proposition
\ref{convmles} concerning the strong topology are still true.}\wbx\end{nota}

In \cite{NeR} the authors proved strong convergence of the orbits of some semigroups.
More than two decades later the result was extended in \cite{Xu} for the almost-orbits.
This extension is straightforward using Proposition \ref{convmles}. An extension to more
general spaces can be found in \cite{GaR}.


Given $v\in L^\infty_{\hbox{\tiny loc}}(0,\infty;X)$, define
$\overline{v}(t)=\frac{1}{t}\int_0^tv(\xi)\ d\xi$. We say $v$ converges strongly (weakly)
{\em in average} if $\overline{v}(t)$ has a strong (weak) limit as $t\to\infty$.
Convergence in average is also inherited by almost-orbits.

\begin{nota}\label{R:translation}{\em Given $v\in L^\infty_{\hbox{\tiny loc}}(0,\infty;X)$, $h,t\geq
0$, define $v_h(t)=v(h+t)$. Since
$\overline{v_h}(t)=\frac{1}{t}\int_0^tv(h+\xi)d\xi=\left(\frac{t+h}{t}\right)\frac{1}{t+h}\int_0^{t+h}v(\eta)d\eta
-\frac{1}{t}\int_0^hv(\xi)d\xi$, if $v$ converges strongly (weakly) in average to $L$,
the same holds for $v_h$, for each $h\ge 0$. }\wbx\end{nota}

\begin{prop}\label{averagemles}
Let $U$ be an \mles. If every orbit of $U$ converges strongly (weakly) in average, so does every almost-orbit.
\end{prop}

\dem Let $u$ be an almost-orbit of $U$. For $p,h\ge 0$ and $t$ sufficiently large, define
$\sigma_h(t,p)=\frac{1}{t}\int_0^tU(p+h+\xi,p)u(p)\ d\xi$ and set $\zeta(p)=\tau-\limty{t}\sigma_0(t,p)$, where $\tau$ stands for either the strong or the weak topology according to the hypothesis. Notice that
\begin{equation}\label{E:average}
[\sigma_0(t,p+h)-\sigma_0(t+h,p)]-\left[\sigma_h(t,p)-\sigma_0(t+h,p)\right]=\left[\sigma_0(t,p+h)-\sigma_h(t,p)\right].
\end{equation}
By virtue of Remark \ref{R:translation},
$\tau-\limty{t}\sigma_h(t,p)=\tau-\limty{t}\sigma_0(t,p)=\tau-\limty{t}\sigma_0(t+h,p)$
for each $h\ge 0$. We let $t\to\infty$ in equation \eqref{E:average} and use the weak lower-semicontinuity of the norm to obtain
$\|\zeta(p+h)-\zeta(p)\|\le \|\sigma_0(t,p+h)-\sigma_h(t,p)\|\le M\ \|u(p+h)-U(p+h,p)u(p)\|,$
which in turn tends to zero as $p\to\infty$ uniformly in $h\ge 0$. As a consequence, $\zeta(p)$ converges strongly to some $\zeta_\infty$ as $p\to\infty$. Finally, for any $p,h\geq 0$ we write
\vskip-20pt\begin{eqnarray*}
\overline{u}(p+h)-\zeta_\infty & = & \frac{1}{p+h}\int_0^pu(\xi)\ d\xi +\left[\frac{h}{p+h}\sigma(h,p)-\zeta(p)\right]+[\zeta(p)-\zeta_\infty]\\
&&+\hskip10pt\frac{1}{p+h}\int_0^{h}[u(p+\xi)-U(p+\xi,p)u(p)]\ d\xi.
\end{eqnarray*}
The second term is bounded by $\sup_{k\ge 0}\|u(p+k)-U(p+k,p)u(p)\|$, which is
independent of $h$ and tends to zero as $p\to \infty$. The last term converges strongly
to zero as $p\to \infty$. Thus, given any $\eps>0$, we can choose $p_\eps$ large enough
so that the second and forth terms are both less than $\eps$. Having fixed $p_\eps $, the
first term converges strongly to zero as $h\to\infty$ while the third term
$\tau$-converges to zero. As a consequence $\overline{u}(t)$ is $\tau$-convergent to
$\zeta_\infty$ as $t\to\infty$. \bx


A function $v\in L^\infty_{\hbox{\tiny loc}}(0,\infty; X)$ is strongly (weakly) {\em
almost-convergent} (in the sense of Lorentz, \cite{Lor}) if there is $y\in X$ such that
$\overline{v_h}(t)$ converges strongly (weakly) to $y$ as $t\to\infty$ uniformly in $h\ge
0$. Almost-convergence implies convergence in average. Conversely, according to Remark
\ref{R:translation}, if $v$ converges in average then $\overline{v_h}$ converges for each
$h\ge 0$, so the uniformity in $h\ge 0$ is what makes the difference. Almost-convergence
is interesting because a trajectory $v(t)$ is convergent if, and only if, it is
almost-convergent and {\em asymptotically regular} (the difference $v(t+h)-v(t)$
converges to zero as $t\to\infty$ for each $h\ge 0$) for the corresponding topology (see
\cite{Lor}). This fact $-$ or method of proof $-$ has been used, for instance, in
\cite{BBR}. The following result shows that almost-convergence of the orbits is also
inherited by the almost-orbits.

\begin{prop}\label{almostmles}
Let $U$ be an \mles. If every bounded orbit of $U$ is strongly (weakly)
almost-convergent, so is every bounded almost-orbit of $U$.
\end{prop}

\dem Define $\sigma_h(t,p)=\frac{1}{t}\int_0^tU(p+h+\xi,p)u(p)\ d\xi $, where $u$ is a bounded almost-orbit of $U$. According to Remark \ref{boundbound}, there exists $p_0\ge 0 $ such that  for all $p\geq p_0$ and $h\geq 0$ we have $\|U(p+h,p)u(p)\|\leq 1+\|u\|_\infty.$ Therefore, by virtue of the hypothesis, for every $p\geq p_0$ there exists $\zeta(p)\in X$ such that for all $h\geq 0$,
$\zeta(p)=\tau-\lim_{t\to\infty}\sigma_h(t,p)$, and the convergence is uniform in $h\ge 0$. Next, we  prove that $\{\zeta(p):p\ge 0\}$ is a Cauchy net. For every $p,h\ge 0$  and $t\geq p+h$ we have $\|\sigma_0(t,p+h)-\sigma_h(t,p)\|\le M\|u(p+h)-U(p+h,p)u(p)\|$.
Let $t\to \infty$ to get $\|\zeta(p+h)-\zeta(p)\|\le M\|u(p+h)-U(p+h,p)u(p)\|,$ which tends to $0$ as $p\to\infty$ uniformly in $h\ge 0$. Hence $\zeta(p)\to\zeta_\infty$ as $p\to\infty$, for some $\zeta_\infty$. For any $p,h,k\geq 0$ we write
\vskip-20pt\begin{eqnarray*}
\overline{u_k}(p+h)-\zeta_\infty & = & \frac{1}{p+h}\int_0^pu(k+\xi)\ d\xi +\left[\frac{h}{p+h}\sigma_k(h,p)-\zeta(p)\right]\\
&&\hskip10pt+ \frac{1}{p+h}\int_0^{h}[u(p+k+\xi)-U(p+k+\xi,p)u(p)]\ d\xi\\
&&\hskip20pt+[\zeta(p)-\zeta_\infty].
\end{eqnarray*}
The first term on the right-hand side is bounded by $p/(p+h)\|u\|_\infty$, independently of $k$. The  second term  is bounded by $\sup_{q\ge
0}\|u(p+q)-U(p+q,p)u(p)\|$, which is independent of $h$ and $k$, and tends to zero as $p\to \infty$. The last term converges strongly to zero as $p\to \infty$. Thus, given any $\eps>0$, we can choose $p_\eps$ large enough so that the second and forth terms are both less than $\eps$. Then, for such  $p_\eps $, the first term converges strongly to zero as $h\to\infty$ while the third term $\tau$-converges to zero, both  uniformly in $k$. As a consequence $\overline{u_k}(t)$ is $\tau$-convergent to $\zeta_\infty$ as $t\to\infty$ uniformly in $k$. \bx

\begin{nota}\label{weaklycomplete}{\em Proposition \ref{almostmles} was proved in \cite{KoM} under
additional assumptions: i) $U$ is an autonomous and strongly continuous \ces, ii) the set
of stationary points is nonempty, and iii) for the weak topology, the space $X$ is
assumed to be {\em weakly complete}, which means that every weak Cauchy net converges
weakly to an element in $X$. The spaces $\luno$, $L^1$ and all reflexive Banach spaces
have this property. It is not the case if $X$ contains $c_0$, though. }\wbx\end{nota}

\begin{nota}{\em Compared with Propositions \ref{convmles} and \ref{averagemles},
the hypotheses and conclusion in Proposition \ref{almostmles} are weaker. According to
Remark \ref{R:all_bounded}, the two formulations are equivalent whenever the ES has
bounded almost-orbits. For Proposition \ref{almostmles} to be useful, one must prove that
the system has at least one bounded almost-orbit. In practice, this step tends to be
useful for proving that the orbits are convergent. In many applications one has to do it
anyway. }\wbx\end{nota}


Let us introduce some general notions of convergence with respect to time-dependent
probability measures in order to unify and summarize the results of the previous section.

Let $\mu$ be a probability measure on $[0,\infty)$. A function $v\in
L^\infty_{\hbox{\tiny loc}}(0,\infty;X)$ is {\em $\mu$-integrable} if the {\em $\mu$-mean
of $v$} on $[0,\infty)$, $\mu(v)=\int_0^{\infty}v(\xi)d\mu(\xi)$ exists. Given a family
$\{\mu_t\}_{t\ge 0}$ of probability measures on $[0,\infty)$, a function $v\in
L^\infty_{\hbox{\tiny loc}}(0,\infty;X)$ is {\em $\{\mu_t\}$-integrable} if $\mu_t(v)$
exists for all $t\geq 0$. We say $v$ {\em converges to $y$ in $\mu_t$-mean} for the
topology $\tau$ if $y=\tau-\lim_{t\to\infty}\mu_t(v)$.

\begin{ejem}\label{cia1}{\em Let $v\in L^\infty_{\hbox{\tiny loc}}(0,\infty;X)$. If
$\mu_t=\delta_t$ is the Dirac mass at $t$, then $\mu_t(v)=v(t)$ and convergence in
$\mu_t$-mean is standard convergence. If
$d\mu_t(\xi)=\hbox{$\frac{1}{t}$}\chi_{[0,t]}(\xi)d\xi$, where $\chi_A$ is the
characteristic function of the set $A$, then
$\mu_t(v)=\hbox{$\frac{1}{t}$}\int_0^tv(\xi)\ d\xi=\overline{v}(t)$ and convergence in
$\mu_t$-mean is convergence in average. }\wbx\end{ejem}

In the rest of this section, $\{\mu_t\}_{t\ge 0}$ is a family of probability measures
such that $\mu_t([0,p])\to 0$ as $t\to \infty$ for all $p\geq0$. Under this assumption,
Propositions \ref{mu-mean} and \ref{mu-mean-uniform} below can be proved by a direct
adaptation of the proofs of Propositions \ref{averagemles} and \ref{almostmles},
respectively. We leave the details to the reader. We restate them as Theorem
\ref{restatement} in terms of almost-equivalent evolution systems.

\begin{prop}\label{mu-mean}
If the $\mu_t$-mean of every orbit of an \mles\ $U$ converges strongly (weakly) as $t\to
\infty$, so does the $\mu_t$-mean of every almost-orbit.
\end{prop}

Given $v\in L^\infty_{\hbox{\tiny loc}}(0,\infty;X)$ and $h\geq 0$, we set
$v_h(t)=v(h+t)$ for $t\geq 0$. If there is $y\in X$ such that
$y=\tau-\lim_{t\to\infty}\mu_t(v_h)=\tau-\lim_{t\to\infty}\int_0^{\infty}v(h+\xi)\
d\mu_t(\xi)$ uniformly in $h\geq0$, for $\tau$ the strong (weak) topology, we say $v$
converges strongly (weakly) to $y$ {\em in $\mu_t$-mean, uniformly with respect to
translations}.

\begin{ejem}{\em \label{e:mu-mean-unif} If $\mu_t$ is the Dirac mass at $t$, then
$\mu_t(v_h)=v(t+h)$ and convergence in $\mu_t$-mean recovers standard convergence, which
is automatically uniform with respect to translations. If
$d\mu_t(\xi)=\hbox{$\frac{1}{t}$}\chi_{[0,t]}(\xi)d\xi$, then
$\mu_t(v_h)=\hbox{$\frac{1}{t}$}\int_0^tv(h+\xi)\ d\xi$. In this case, convergence in
$\mu_t$-mean uniformly with respect to translations is almost-convergence.
}\wbx\end{ejem}

\begin{prop}\label{mu-mean-uniform}
If every bounded orbit of an \mles\ $U$ converges strongly (weakly) in $\mu_t$-mean
uniformly with respect to translations, so does every bounded almost-orbit.
\end{prop}

\begin{thm}\label{restatement}
Let $U$ and $V$ be two \mles\ which are AAE. If every orbit of $U$ converges strongly
(weakly) in $\mu_t$-mean, so does every orbit of $V$. If the convergence is uniform with
respect to translations for all bounded orbits of $U$, the same holds for the bounded
orbits of $V$.
\end{thm}


\section{Further results on Lipschitz evolution systems}

If $U$ is an ES on $C$, let $SP(U):=\{x\in C | U(t,s)x=x, \forall t\geq s\}$ be its
(possibly empty) set of {\em stationary points}. Similarly, denote by $\ASP(U)$ its set
of {\em almost-stationary points}: $\ASP(U)=\{x\in C | \lim_{t\to\infty} \sup_{h\ge
0}\|U(t+h,t)x-x\|=0\}$. Clearly $SP(U) \subseteq \ASP(U)$ and if $U$ is autonomous, then
 $\ASP(U)=SP(U)$. This is not the case in general even for a
\ces\ (take $U(t,s)x=x+e^{-s}-e^{-t}$, $x\in \R$, for which $\ASP(U)=\R$ and
$SP(U)=\emptyset$).

\begin{nota}{\em A nice characterization of $SP(T)$ is given in \cite{Suz} when $T$ is an autonomous
\ces\ on a weakly compact subset of a Banach space with the Opial property: $z\in SP(U)$
if, and only if, there is $t_n\to\infty$ such that
$w$-$\limn\frac{1}{t_n}\int_0^{t_n}T(s)zdz=z$. A similar result using strong limits is
given in \cite{SuT}. A challenging task is to extend this characterization to
nonautonomous ES. }\wbx\end{nota}

If $x^*\in\ASP(U)$, the constant function $u(t)\equiv x^*$ is a bounded almost-orbit of
$U$. According to Remark \ref{R:all_bounded}, if $\ASP(U)\neq\emptyset$ every
almost-orbit of $U$ is bounded. We now turn our attention to closedness and convexity of
$\ASP(U)$.

\begin{lem} Suppose $C$ is closed for the strong topology.
If $U$ is an \mles\ on C then $\ASP(U)$ is closed for the strong topology.
\end{lem}

\dem Let $\{x_n\}\in\ASP(U)$ converge to $x$. Since $x\in C$, $\|U(t+h,t)x-x\|\le
(M+1)\|x_n-x\|+\|U(t+h,t)x_n-x_n\|$. Hence $\limsup\limits_{t\to\infty}\sup\limits_{h\ge
0}\|U(t+h,t)x-x\|\le (M+1)\|x_n-x\|.$ Letting $n\to \infty$, we conclude that
$x\in\ASP(U)$.\bx

By Remark \ref{R:limits}, if $U$ is an \mles, $\ASP(U)$ contains all the limits of the
strongly convergent almost-orbits if there are any. For weak limits we have the
following:

\begin{cor}
Let $C$ be strongly closed. If the weak limits of orbits of an \mles\ $U$ lie in
$\ASP(U)$, the same holds for weak limits of almost-orbits of $U$.
\end{cor}

\dem With the notation introduced in the proof of Proposition \ref{convmles}, as $t\to\infty$, $u(t)$ converges to
$\zeta_\infty$, which is the strong limit of weak limits of orbits of $U$.\bx

\begin{nota}{\em If $U$ is an autonomous \mles\ whose orbits converge weakly to points in $SP(U)$,
the weak limits of almost-orbits are also in $SP(U)$. }\wbx\end{nota}

Let $K\subset X$ be nonempty and convex. Let $\gamma:\R_+\to\R_+$ be continuous and
strictly increasing with $\gamma(0)=0$. A function $F:K\to X$ is {\em of type $\gamma$}
if $\gamma(\|F(\lambda x+(1-\lambda)y)-\lambda F(x)-(1-\lambda)F(y)\|)\le
\|x-y\|-\|F(x)-F(y)\|$ for all $x,y\in K$ and $\lambda\in [0,1]$.


\begin{prop}\label{P:convex}
Let $U$ be a \ces\ on a convex set $C$ and suppose that there exists $\gamma$ such that
for each $t\ge s$, $U(t,s)$ is of type $\gamma$ on a convex set $K$ containing $\ASP(U)$.
Then $\ASP(U)$ is convex.
\end{prop}

\dem Take $x_1,x_2\in\ASP(U)$ and define $\psi_i(t)=\sup_{h\ge 0}\|U(t+h,t)x_i-x_i\|$ for $i=1,2$. Now take
$\lambda\in(0,1)$ and set $z=\lambda x_1+(1-\lambda)x_2$. We have
\vskip-20pt\begin{eqnarray*}
\hskip-10pt\|U(t+h,t)z-z\| & \le & \|U(t+h,t)z-\lambda U(t+h,t)x_1-(1-\lambda)U(t+h,t)x_2\|\\
&& \hskip10pt+\lambda\psi_1(t)+(1-\lambda)\psi_2(t)\\
& \le & \gamma^{-1}\left(\|x_1-x_2\|-\|U(t+h,t)x_1-U(t+h,t)x_2\|\right)\\
& & \hskip10pt+\lambda\psi_1(t)+(1-\lambda)\psi_2(t)\\
& \le & \gamma^{-1}\left(\psi_1(t)+\psi_2(t)\right)+\lambda\psi_1(t)+(1-\lambda)\psi_2(t)
\end{eqnarray*}
\vskip-20pt Letting $t\to\infty$ we get the result.\bx

If $X$ is uniformly convex and $K$ is bounded and convex, there is $\gamma$ such that
every nonexpansive function $F:K\to X$ is of type $\gamma$ (see \cite[Lemma 1.1]{Bru4}).
If $U$ is a \ces\ on a convex set $C$ then for every bounded $A\subset C$, there is
$\gamma$ such that $U(t,s)$ is of type $\gamma$ on ${\rm co}(A)$ for each $t\ge s$. We
deduce the following:

\begin{cor}
Let $U$ be a \ces\ on a convex subset $C$ of a uniformly convex Banach space $X$. Then $\ASP(U)$ is convex.
\end{cor}

\dem Let $x_1,x_2\in\ASP(U)$. The line segment $K$ joining $x_1$ and $x_2$ is bounded and
convex, so there is a function $\gamma$ such that the restriction of $U(t,s)$ to $K$ is
of type $\gamma$ for all $t\ge s\ge 0$. The rest follows as in Proposition
\ref{P:convex}.\bx

\begin{nota}{\em Proposition \ref{P:convex} is valid if we replace $\ASP(U)$ with the set of all
almost-orbits. The proof is quite similar and uses Proposition \ref{prop1} iv)  for
concluding. As a consequence, if $X$ is uniformly convex and  $U$ is a \ces\ on a bounded
convex set $C$, then the set of all almost-orbits is convex. }\wbx\end{nota}



Let $\{M(t,s)\}_{t\ge s\ge 0}$ be a family of positive numbers satisfying
$\limty{t}M(t,s)=0$ for each $s$. A {\em strongly contracting evolution system} (\sces)
on $C$ is an ES $U$ such that $\|U(t,s)x-U(t,s)y\|\le M(t,s)\|x-y\|$ for all $x,y\in C$
and $t\ge s\ge 0$.


\begin{prop}\label{propsces}
Let $U$ be a \sces. We have the following:
\begin{itemize}
    \item [i)] If $u_1$ and $u_2$ are almost-orbits of $U$ then $\limty{t}\|u_1(t)-u_2(t)\|=0$;
    \item [ii)] The set $\ASP(U)$ has at most one element; and
    \item [iii)] If $\ASP(U)\neq\emptyset$, then every almost-orbit of $U$ converges strongly to the unique $x^*\in\ASP(U)$.
\end{itemize}
\end{prop}

\dem We have $\|u_1(t+s)-u_2(t+s)\| \le \psi_1(t)+\psi_2(t)+M(t+s,t)\|u_1(t)-u_2(t)\|$ as in part $i)$ of Proposition \ref{prop1}.
But $\limsup\limits_{s\to\infty}\|u_1(s)-u_2(s)\|\le \psi_1(t)+\psi_2(t)$, so $\limty{s}\|u_1(s)-u_2(s)\|=0$.
Parts ii) and iii) are a trivial consequence.\bx

\section{Asymptotic equivalence without the Lipschitz condition}

Not all the results in the previous section are true without the Lipschitz assumption on
the ES. For instance, having a bounded almost-orbit does not imply that all the
almost-orbits are bounded (take $U(t,s)x=e^{(t-s)}x$). We discuss on some properties of
the orbits that do hold for the almost-orbits. The first work that contains equivalence
results for non-contracting ES seems to be \cite{Miy}, where they study strongly
continuous semigroups which are ``asymptotically nonexpansive in the intermediate sense".
We shall not go into the details but just mention that the existing results on asymptotic
equivalence require additional (and strong!) regularity with respect both to time and
space.


Let $\{\mu_t\}_{t\geq 0}$ be a family of probability measures on $[0,\infty)$.

\hangindent=20mm\hangafter=-3\vspace{-3mm}
\hskip-20mm Hypothesis $\mathbf{H}$: For each  $\{\mu_t\}$-integrable $g$ with $\limty{t}\int_0^\infty g(\xi)\ d\mu_t(\xi)=L$, each
$\eps>0$ and $K>0$ there exists $T>0$ such that for all $t\ge T$ one
has $\left\|\frac{}{}\int_0^\infty g(\xi)\ d\mu_t(\xi+K)-L\right\|<\eps.$

The families described in Example \ref{cia1} do satisfy Hypothesis {\bf H}: This is
trivial if $\mu_t$ is the Dirac mass at $t$. If
$d\mu_t(\xi)=\hbox{$\frac{1}{t}$}\chi_{[0,t]}(\xi)$, then for $t$ large enough
$\int_0^\infty
g(\xi)d\mu_t(\xi+K)=\left(\frac{t-K}{t}\right)\frac{1}{t-K}\int_0^{t-K}g(\xi)d\xi,$ which
tends to $L$ as $t\to\infty$. The fact that $\limty{t}\mu_t(B)=0$ for each bounded set
$B$ does not imply that Hypothesis $\mathbf{H}$ will hold:

\begin{ejem}{\em  Define $n(\xi)=\sum_{k\ge 0}\chi_{[2k,2k+1)}(\xi)$ and $\hat n(\xi)=n(\xi+1)$ so
that $n^2\equiv n$ and $n\hat n\equiv 0$. Let
$d\mu_t(\xi)=\alpha^{-1}(t)n(\xi)\chi_{[0,t]}(\xi)d\xi$, where
$\alpha(t)=\int_0^tn(\xi)d\xi$. Then $\mu_t(B)\to 0$ for every bounded set $B$ (this is
obvious) but does not fulfill Hypothesis {\bf H}. To see this, simply notice that
$\int_0^\infty n(\xi)d\mu_t(\xi)=1$ while $\int_0^\infty
n(\xi)d\mu_t(\xi+1)=\alpha^{-1}\int_1^{t-1} \hat n(\xi)n(\xi)d\xi=0$ for all
$t$.}\wbx\end{ejem}

\begin{thm}\label{teocia1}
Let $U$ be an ES and let $\{\mu_t\}$ satisfy Hypothesis {\bf H}. If each orbit converges
strongly in $\mu_t$-mean, so does every $\{\mu_t\}$-integrable almost-orbit.
\end{thm}

\dem Suppose $u$ is a $\{\mu_t\}$-integrable almost-orbit of $U$ and let $\eps>0$. Choose
$S>0$ such that $\sup_{h\ge 0}\|u(t+h)-U(t+h,t)u(t)\| <\eps/6$ for all $t\ge S$. Define
$\zeta(S)=\limty{t}\int_0^\infty U(S+\xi,S)u(S)\ d\mu_t(\xi)$. By hypothesis, there is
$T_1$ such that $\left\|\zeta(S)-\int_0^\infty U(S+\xi,S)u(S)\
d\mu_t(\xi)\right\|<\eps/6$ for all $t\ge T_1$. We have\\
$\left\|\mu_t(u)-\zeta(S)\right\| \le \int_0^{S}\left\|u(\xi)\right\|\
d\mu_t(\xi)+\int_{S}^{\infty}\|u(\xi)-U(\xi,S)u(S)\|\ d\mu_t(\xi)$\\
$\frac{}{}$\hfill $+ \left\|\zeta(S)-\int_{0}^{\infty}U(S+\xi,S)u(S)\
d\mu_t(\xi+S)\right\|.$\\
For the first term, since $\limty{t}\mu_t([0,S])= 0$, we can take $T_2$ such that
$\mu_t([0,S])<\eps/6C$ for all $t\ge T_2$, where $C=\sup_{0\le\xi\le S}\|u(\xi)\|$. The
second term is less than $\eps/6$. By Hypothesis {\bf H} there is $T_3$ such that the
last term is less than $\eps/6$ whenever $t\ge T_3$. Hence if $t\ge T=\max\{\ T_1,\ T_2,\
T_3\ \}$, we have $\|\mu_t(u)-\zeta(S)\|<\eps/2$ for all $h\ge 0$. We have found $T>0$
such that $\|\mu_t(u)-\mu_s(u)\|<\eps$ for all $t,s\ge T$ and therefore $\mu_t(u)$
converges to some $y$ as $t\to\infty$.\bx

If $\mu_t=\delta_t$ the argument above gives an alternative proof of the assertion for
the strong topology in Proposition \ref{convmles} without the Lipschitz assumption.
Moreover, this proof is simpler because we do not perform the intermediate step of
proving the convergence of $\zeta(p)$. This argument fails when dealing with the weak
topology. To overcome this problem, consider the following hypothesis, also satisfied by
the families described in Example \ref{cia1}:

\hangindent=20mm\hangafter=-3\vspace{-3mm}
\hskip-20mm Hypothesis {\bf w-H}: For each $\{\mu_t\}$-integrable $g$ with w-$\limty{t}\int_0^\infty g(\xi)\ d\mu_t(\xi)=L$,
 each $\eps>0$, $K>0$ and $f\in X^*$ there exists $T>0$ such that for all
$t\ge T$ one has $\left|\langle\frac{}{}\int_0^\infty g(\xi)\ d\mu_t(\xi+K)-L,f\ \rangle\right|<\eps.$

Under Hypothesis {\bf w-H} the argument above shows that if the orbits of $U$ converge
weakly in $\mu_t$-mean, then $\mu_t(u)$ has the Cauchy property for the weak topology
whenever $u$ is an almost-orbit of $U$ (i.e.
$\limty{t,s}\langle\mu_t(u)-\mu_s(u),\phi\rangle=0$ for each $\phi\in X^*$). If $X$ is
weakly complete (see Remark \ref{weaklycomplete}) the net $\{\mu_t(u)\}$ converges
weakly. This is a version of Theorem \ref{teocia1} for the weak topology.

Recall that a $\{\mu_t\}$-integrable function $v$ $\tau$-converges to $y\in X$ in
$\mu_t$-mean, uniformly with respect to translations if $\mu_t(v_h)$ $\tau$-converges to
$y$ as $t\to\infty$ uniformly in $h\ge 0$. This notion includes standard convergence and
almost-convergence for the families of measures in Example \ref{cia1}. The uniformity in
$h\ge 0$ requires a slightly stronger assumption on $\{\mu_t\}$ (that still hold for the
families mentioned above) in order to prove the equivalence results:

\hangindent=20mm\hangafter=-3\vspace{-3mm}
\hskip-20mm Hypothesis $\mathbf{H_u}$: For each $\{\mu_t\}$-integrable $g$ with $\limty{t}\int_0^\infty g(\xi)\ d\mu_t(\xi)=L$,
each $\eps>0$ and $K>0$ there exists $T>0$ such that for all $t\ge T$
and $k\in[0,K]$ one has $\left\|\frac{}{}\int_0^\infty g(\xi)\ d\mu_t(\xi+k)-L\right\|<\eps.$

\begin{thm}
Let $U$ be an ES and assume $\{\mu_t\}$ satisfies Hypothesis $\mathbf{H_u}$. If $U(t,s)x$
converges strong\-ly in $\mu_t$-mean, uniformly with respect to translations for all $x$
and $s$, then so does every $\{\mu_t\}$-integrable almost-orbit.
\end{thm}

\dem Suppose $u$ is a $\{\mu_t\}$-integrable almost-orbit of $U$ and let $\eps>0$. Choose
$S>0$ such that $\sup\limits_{h\ge 0}\|u(t+h)-U(t+h,t)u(t)\| <\eps/6$ for all $t\ge S$.
Define $\zeta(S)=\limty{t}\int_0^\infty U(S+\xi,S)u(S)\ d\mu_t(\xi)$. By hypothesis,
there is $T_1$ such that $\left\|\zeta(S)-\int_0^\infty U(S+h+\xi,S)u(S)\
d\mu_t(\xi)\right\|<\eps/6$ for all $t\ge T_1$ and $h\ge 0$ (the convergence is uniform
in $h\ge 0$). We divide the rest of the proof in two parts:

\underline{$0\le h\le S$}: As in the proof of Theorem \ref{teocia1} we have\\
$\left\|\mu_t(u_h)-\zeta(S)\right\|\le \int\limits_0^{S-h}\!\!\!\left\|u(h+\xi)\right\|
d\mu_t(\xi)
+\int\limits_{S-h}^{\infty}\!\!\!\|u(h+\xi)-U(h+\xi,S)u(S)\| d\mu_t(\xi)$\\
$\frac{}{}$\hfill $+ \left\|\zeta(S)-\int_{0}^{\infty}U(S+\xi,S)u(S)\
d\mu_t(\xi+(S-h))\right\|$.\\
For the first term, since $\mu_t([0,S])\to 0$ as $t\to\infty$, we can take $T_2$ such
that $\mu_t([0,S])<\eps/6C$ for all $t\ge T_2$, where $C=\sup_{0\le\xi\le S}\|u(\xi)\|$.
The second term is always less than $\eps/6$. Finally, use Hypothesis $\mathbf{H_u}$ to
find $T_3$ such that the last term is less than $\eps/6$ whenever $t\ge T_3$. Hence if
$t\ge T=\max\{\ T_1,\ T_2,\ T_3\ \}$, we have $\|\mu_t(u_h)-\zeta(S)\|<\eps/2$ for all
$h\ge 0$.\\
$\underline{h\ge S}:\ \left\|\mu_t(u_h)-\zeta(S)\right\|  \le
\int\limits_0^\infty\|u(h+\xi)-U(h+\xi,S)u(S)\|\ d\mu_t(\xi)$\\
$\frac{}{}$\hfill $+ \left\|\zeta(S)-\int_0^\infty U(h+\xi,S)u(S)\
d\mu_t(\xi)\right\|$,\\
whenever $t\ge T_1$. Each term is less than $\eps/6$, so
$\|\mu_t(u_h)-\zeta(S)\|<\eps/3<\eps/2$ for all $t\ge T_1$ and $h\ge S$. Finally,
$\|\mu_t(u_h)-\zeta(S)\|<\eps/2$ for all $t\ge T$ and $h\ge 0$. This implies
$\|\mu_t(u_h)-\mu_s(u_k)\|<\eps$ for all $t,s\ge T$ and $h,k\ge 0$ and so $u$ is strongly
convergent in $\mu_t$-mean, uniformly with respect to translations. \bx


We obtain the corresponding result for the weak topology (essentially with the same proof
if $X$ is weakly complete) if we replace Hypothesis $\mathbf{H_u}$ by:

\hangindent=20mm\hangafter=-3\vspace{-3mm}
\hskip-20mm Hypothesis {\bf w-}$\mathbf{H_u}$: For each $\{\mu_t\}$-integrable $g$ with w-$\limty{t}\int_0^\infty g(\xi)\ d\mu_t(\xi)=L$,
each $\eps>0$, $K>0$ and $f\in X^*$ there exists $T>0$ such that for all
$t\ge T$ and $k\in[0,K]$ one has $\left|\langle\frac{}{}\int_0^\infty g(\xi)\ d\mu_t(\xi+k)-L,f\ \rangle\right|<\eps.$


\section{A couple of additional results}

For an autonomous \ces, every orbit is uniformly continuous. According to \cite{KoM}, so
is every continuous almost-orbit. Their proof uses the contracting property, which we
show to be unnecessary. The key lies on the time-dependence.

\begin{prop}
Let $U$ be an evolution system. If every orbit is uniformly continuous, so is every continuous almost-orbit.
\end{prop}

\dem Let $u$ be a continuous almost-orbit of $U$ and $\eps>0$. First, take $T>0$ such that
$\|u(\tau)-U(\tau,T)u(T)\|<\eps/3$ for all $\tau\ge T$. Since $u$ is continuous, it is uniformly continuous on
$[0,T+1]$. Hence there is $\delta_1>0$ such that for every $t,s\in[0,T+1]$ satisfying $|t-s|<\delta_1$ one has
$\|u(t)-u(s)\|<\eps$. Now consider the function $\tau\mapsto U(\tau,T)u(T)$ defined for $\tau\ge T$. By hypothesis it
is uniformly continuous, so there is $\delta_2>0$ such that for every $t,s\ge T$ such that $|t-s|<\delta_2$ one has
$\|v(t,T)u(T)-U(s,T)u(T)\|<\eps/3$. Therefore, if $t,s\ge T$ and $|t-s|<\delta_2$, the quantity $\|u(t)-u(s)\|$ is bounded by
$\|u(t)-U(t,T)u(T)\|+\|U(s,T)u(T)-u(s)\|+\|U(t,T)u(T)-U(s,T)u(T)\|.$ For $t,s\in[0,\infty)$ with $|t-s|<\min\{\delta_1,\delta_2,1\}$ we have $\|u(t)-u(s)\|<\eps$. \bx

The $\omega$-limits of almost-orbits have some kind of invariance under $U$:

\begin{prop}
Let $u$ be an almost-orbit of an ES $U$ and suppose $\{s_n\}$ is strictly increasing with
$\limn s_n=\infty$ and $\tau$-$\limn u(s_n)=x^*$.
\begin{itemize}
    \item [i)] There exist a sequence $\{t_n\}$ of positive numbers and a sequence $\{x_n\}$ in $C$ such that
        $\tau$-$\limn x_n=x^*$, $t_n>s_n$ and $\tau$-$\limn U(t_n,s_n)x_n=x^*$. The sequence
        $\{t_n\}$ can be chosen such that $\limn(t_n-s_n)=\infty$.
    \item [ii)] If $U$ is an \mles\, and $\tau$ is the strong topology, then $\limn U(t_n,s_n)x^*=x^*$.
\end{itemize}
\end{prop}

\dem For the first part, let $\f:\N\to\N$ be any positive function
and set $h_n=s_{n+\f(n)}-s_n$. Write $t_n=s_n+h_n$ and $x_n=u(s_n)$.
Since $u$ is an almost orbit of $V$, for every $\eps>0$ there is
$N\ge0$ such that $\|u(s_n+h_n)-U(s_n+h_n,s_n)u(s_n)\|<\eps$ for all
$n\ge N$. Therefore $U(t_n,s_n)x_n-x^*=U(s_n+h_n,s_n)u(s_n)-u(s_n+h_n)+u(s_{n+\f(n)})-x^*$
tends to zero for the topology $\tau$. Clearly $\f$ can be chosen so that $s_{n+\f(n)}-s_n$ tends to $\infty$ as $n\to\infty$. For the second part, just notice that
$\|U(t_n,s_n)x^*-x^*\| \le  M\|u(s_n)-x^*\| + \|U(t_n,s_n)u(s_n)-u(t_n)\| + \|u(t_n)-x^*\|,$
which tends to zero as $n\to\infty$. \bx

\section{Concluding remarks}

The tools developed here are potentially useful in different scenarios, namely: in
general {\em asymptotic analysis}, information on the asymptotic behavior of a system can
be derived from the study of one that is AAE (as in \cite{Pas} and \cite{Gul}). In {\em
numerical analysis}, to determine whether a discretization has the same asymptotic
properties as the continuous-time model. For instance, it would be possible to know {\em
a priori} if one must take averages in order to approximate the solution of a problem. In
{\em perturbations theory}, to know how much a system can be perturbed without changing
its asymptotic behavior. This could help predict or control the effect of errors and
noises. For {\em ill-posed problems}, to get an idea of what kind of perturbations can
force a system to converge when it does not. For example, in some optimization problems
it is known that a viscosity term can force a nonconverging system to converge (see
\cite{AtC} or \cite{ACz}). Several applications in optimization, fixed-point theory,
games theory and parabolic equations will be presented in a forthcoming paper in
preparation.

\end{document}